\documentclass[12pt]{amsart}

\usepackage{enumerate}
\usepackage{amsmath,xspace,amssymb,mathrsfs}
\usepackage{color}
\definecolor{darkgreen}{rgb}{0,0.45,0} 
\usepackage[pagebackref,colorlinks,citecolor=darkgreen,linkcolor=darkgreen]{hyperref}

\input xy
\xyoption{all}
\xyoption{2cell}
\UseAllTwocells
\CompileMatrices

\renewcommand{\phi}{\varphi}

\newcommand{\C}{{\ensuremath{\mathscr C}}\xspace}
\newcommand{\D}{{\ensuremath{\mathscr D}}\xspace}

\newcommand{\E}{{\ensuremath{\mathscr E}}\xspace}
\newcommand{\F}{{\ensuremath{\mathscr F}}\xspace}
\newcommand{\J}{{\ensuremath{\mathscr J}}\xspace}
\newcommand{\K}{{\ensuremath{\mathscr K}}\xspace}

\newcommand{\Q}{{\ensuremath{\mathscr Q}}\xspace}
\newcommand{\R}{{\ensuremath{\mathscr R}}\xspace}
\renewcommand{\S}{{\ensuremath{\mathscr S}}\xspace}
\newcommand{\op}{\ensuremath{{}^{\textrm{op}}}\xspace}
\newcommand{\Set}{\textnormal{\bf Set}\xspace}
\newcommand{\RGph}{\textnormal{\bf RGph}\xspace}
\newcommand{\Preord}{\textnormal{\bf Preord}\xspace}

\DeclareMathOperator{\Sh}{Sh}
\DeclareMathOperator{\Sep}{Sep}
\DeclareMathOperator{\ev}{ev}
\DeclareMathOperator{\Cs}{Cs}

\def\two{\mathbf 2}

\def\x{\times}

  \newtheorem{proposition}[subsection]{Proposition}
  \newtheorem{lemma}[subsection]{Lemma}
  
  \newtheorem{corollary}[subsection]{Corollary}
  \newtheorem{theorem}[subsection]{Theorem}

  \theoremstyle{definition}
  
  \newtheorem{example}[subsection]{Example}
  
  \newtheorem{remark}[subsection]{Remark}

\def\clap#1{\hbox to 0pt{\hss#1\hss}}
\def\mathllap{\mathpalette\mathllapinternal} \def\mathrlap{\mathpalette\mathrlapinternal} 
\def\mathllapinternal#1#2{%
 \llap{$\mathsurround=0pt#1{#2}$}}
\def\mathrlapinternal#1#2{%
\rlap{$\mathsurround=0pt#1{#2}$}}

\begin{document}

\title{Grothendieck quasitoposes}

\author{Richard Garner}
\address{Department of Computing, Macquarie University, NSW 2109 Australia.}
\email{richard.garner@mq.edu.au}

\author{Stephen Lack}
\address{Department of Mathematics, Macquarie University, NSW 2109 Australia.}
\email{steve.lack@mq.edu.au}

\date{November 2011}


\begin{abstract}
A full reflective subcategory \E of a presheaf category $[\C\op,\Set]$ 
is the category  of sheaves for a topology $j$ on \C if and
only if the reflection from $[\C\op,\Set]$ into \E preserves finite limits.
Such an \E is then called a Grothendieck topos. More generally, one can
consider two topologies, $j\subseteq k$, and
the category of sheaves for $j$ which are also separated for $k$. 
The categories \E of this form for some \C, $j$, and $k$ are the
Grothendieck quasitoposes of the title, previously studied by Borceux and 
Pedicchio, and include many examples of categories of spaces. They also
include the category of concrete sheaves for a concrete site. 
We show that a full reflective subcategory \E of $[\C\op,\Set]$ arises 
in this way for some $j$ and $k$ if and only if
the reflection preserves 
monomorphisms as well as pullbacks over elements of \E. 
More generally, for any quasitopos \S, 
we define a subquasitopos of \S to be
a full reflective subcategory of \S for which the reflection preserves 
monomorphisms as well as pullbacks
over objects in the subcategory, and we characterize such 
subquasitoposes in terms of universal closure operators. 
\end{abstract}

\maketitle

\section{Introduction}

A Grothendieck topos is a category of the form $\Sh(\C,j)$ for a small
category \C and a (Grothendieck) topology $j$ on \C. These categories
have been of fundamental importance in geometry, logic, and other 
areas. Such categories were characterized by Giraud as the cocomplete
categories with a generator, satisfying various exactness conditions 
expressing compatibility between limits and colimits. 

The category 
$\Sh(\C,j)$ is a full subcategory of the presheaf category $[\C\op,\Set]$,
and the inclusion has a finite-limit-preserving left adjoint. This in fact
leads to another characterization of Grothendieck toposes. A full 
subcategory is said to be {\em reflective} if the inclusion has a left 
adjoint, and is said to be a {\em localization} if moreover this left adjoint
preserves finite limits. A category is a Grothendieck topos if and only
if it is a localization of some presheaf category $[\C\op,\Set]$ on 
a small category \C. (We shall henceforth only consider presheaves on 
small categories.)

Elementary toposes, introduced by Lawvere and Tierney, generalize 
Grothendieck toposes; the non-elementary conditions of cocompleteness
and a generator in the Giraud characterization are replaced by the requirement
that certain functors, which the Giraud conditions guarantee are continuous,
must in fact be representable.
Yet another characterization of the Grothendieck toposes is as the 
elementary toposes which are locally presentable, in the
sense of \cite{Gabriel-Ulmer}. We cite the encyclopaedic 
\cite{elephant1,elephant2} as a general reference for topos-theoretic
matters.

A {\em quasitopos} \cite{Penon-quasitopos} is a generalization of the 
notion of elementary topos. The main difference is that a quasitopos
need not be balanced: this means that in a quasitopos a morphism may
be both an epimorphism and a monomorphism without being invertible. 
Rather than a classifier for all subobjects,
there is only a classifier for {\em strong} subobjects (see Section~\ref{sect:prelims} below). The definition, then,
of a quasitopos is a category \E with finite limits and colimits, for which
\E and each slice category $\E/E$ of \E is cartesian closed, and 
which has a classifier for strong subobjects. A simple example of a 
quasitopos which is not a topos is a Heyting algebra, seen as a category
by taking the objects to be the elements of the Heyting algebra, with a 
unique arrow from $x$ to $y$ just when $x\le y$. Other examples
include the category of 
convergence spaces in the sense of  Choquet, 
or various categories of differentiable spaces, studied by Chen. 
See \cite{elephant1} once again for generalities about quasitoposes,
and \cite{BaezHoffnung} for the examples involving differentiable 
spaces.

The notion of Grothendieck quasitopos was introduced in 
\cite{BorceuxPedicchio}. Once again, there are various possible
characterizations:

\begin{enumerate}[($i$)]
\item the locally presentable quasitoposes;
\item the locally presentable categories which are locally cartesian 
closed and in which every strong equivalence relation is the kernel
pair of its coequalizer;
\item the categories of the form $\Sep(k)\cap\Sh(j)$ for topologies
$j$ and $k$ on a small category \C, with $j\subseteq k$. 
\end{enumerate}
In ($ii$), an equivalence relation in a category \E is a pair 
$d,c\colon R\rightrightarrows A$ inducing an equivalence relation $\E(X,R)$ on each 
hom-set $\E(X,A)$; it is said to be {\em strong} if  the induced map 
$R\to A\x A$ is a 
strong monomorphism. In ($iii$), we write $\Sh(j)$ for the sheaves
for $j$, and $\Sep(k)$ for the category of separated objects for $k$;
these are defined like sheaves, except that in the sheaf condition we 
ask only for the uniqueness, not the existence, of the gluing. A category \C 
equipped with topologies $j$ and $k$ with $j\subseteq k$ is called 
a {\em bisite} in \cite{elephant2}, and a presheaf on \C which is a 
$j$-sheaf and $k$-separated is then said to be $(j,k)$-{\em biseparated}. 

A special case is where \C has a terminal object and the representable
functor $\C(1,-)$ is faithful, and $k$ is the topology generated by the 
covering families consisting, for each $C\in\C$, of the totality of maps 
$1\to C$.  If $j$ is any subcanonical topology 
contained in $k$, then $(\C,j)$ is a concrete site in the sense
of \cite{BaezHoffnung} (see also 
\cite{Dubuc-ConcreteQuasitopoi,DubucEspanol}) 
for which the concrete sheaves are 
exactly the $(j,k)$-biseparated presheaves.

In the case of Grothendieck toposes, a full reflective subcategory 
of a presheaf category $[\C\op,\Set]$ has the form $\Sh(j)$ for 
some (necessarily unique) topology $j$ if and only if the reflection 
preserves finite limits. The lack of a corresponding result for
Grothendieck quasitoposes is a noticeable gap in the existing theory,
and it is precisely this gap which we aim to fill.

It is well-known that the reflection from $[\C\op,\Set]$ to $\Sep(k)\cap\Sh(j)$ preserves finite products and monomorphisms.
In Example~\ref{ex:main} below, we show that this does not suffice to
characterize such reflections, using
the reflection of directed graphs into preorders as a counterexample. 
We provide a remedy for this in Theorem~\ref{thm:main}, where we
show  that a reflection $L\colon [\C\op,\Set]\to\E$
has this form for topologies $j$ and $k$ if and only if $L$ preserves finite 
products and monomorphisms and is also {\em semi-left-exact}, in the sense of \cite{CassidyHebertKelly}. Alternatively, such $L$ can be characterized
as those which preserve monomorphisms and have {\em stable units},
again in the sense of \cite{CassidyHebertKelly}. The stable unit condition
can most easily be stated by saying that $L$ preserves all pullbacks
in $[\C\op,\Set]$ over objects in the subcategory \E. For each object
$X\in\E$, the slice category $\E/X$ is a full subcategory of 
$[\C\op,\Set]/X$, with a reflection $L_X\colon[\C\op,\Set]/X\to\E/X$ 
given on objects by the action of $L$ on a morphism into $X$; the 
condition that $L$ preserve all pullbacks over objects of \E is 
equivalently the condition that each $L_X$ preserve finite products. 

Since a {\em subtopos} of a topos \S is by definition a full reflective
subcategory of \S for which the reflection preserves finite limits,
the Grothen\-dieck toposes are precisely the subtoposes
of presheaf toposes. Sub\-toposes of an arbitrary topos can be 
characterized in terms 
of Lawvere-Tierney topologies; more importantly for our purposes,
they can be characterized in terms of universal closure operators.

By analogy with this case, we define a  {\em subquasitopos} of a 
quasitopos \S to be a 
full reflective subcategory of \S for which the reflection preserves
monomorphisms and has stable units. Thus a Grothendieck 
quasitopos is precisely a subquasitopos of a presheaf topos. 
We also give a characterization of subquasitoposes of an 
arbitrary quasitopos \S, using universal closure operators. 

We begin, in the following section, by recalling a few basic notions
that will be used in the rest of the paper; then in Section~\ref{sect:conditions}
we study various weakenings of finite-limit-preservation for a 
reflection, and the relationships between these. In Section~\ref{sect:qtop}
we study conditions under which reflective subcategories of quasitoposes
are quasitoposes.
In Section~\ref{sect:monics} we characterize subquasitoposes of
a general quasitopos, before turning, in Section~\ref{sect:Grothendieck}, to subquasitoposes of presheaf toposes and their relationship
with Grothendieck quasitoposes.

\subsection*{Acknowledgements}

We are grateful to the anonymous referee for several 
helpful comments on a preliminary version of the paper, in 
particular for suggesting the formulation of Theorem~\ref{thm:k},
which is more precise than the treatment given in an earlier 
version of the paper. 

This research was supported under the Australian Research
Council's {\em Discovery Projects} funding scheme, 
project numbers DP110102360 (Garner) and DP1094883 (Lack). 

\section{Preliminaries}\label{sect:prelims}

We recall a few basic notions that will be used in the rest of the paper. 

A monomorphism $m\colon X\to Y$ is said to be {\em strong} if for all
commutative diagrams
$$\xymatrix{
X' \ar[r]^{e} \ar[d] & Y' \ar[d] \\ X \ar[r]_{m} & Y }$$
with $e$ an epimorphism, there is a unique map $Y'\to X$ making
the two triangles commute. {\em Strong epimorphisms} are defined dually. 
A strong epimorphism which is also a monomorphism is invertible, and
dually a strong monomorphism which is also an epimorphism is invertible.

A {\em weak subobject classifier} is a morphism $t\colon 1\to\Omega$
with the property that for any strong monomorphism $m\colon X\to Y$
there is a unique map $f\colon Y\to \Omega$ for which the diagram
$$\xymatrix{
X \ar[r]^{m} \ar[d] & Y \ar[d]^{f} \\ 1 \ar[r]_{t} & \Omega }$$
is a pullback.

A category with finite limits is said to be {\em regular} if every morphism 
factorizes as a strong epimorphism followed by a monomorphism, and 
if moreover any pullback of a strong epimorphism is again a strong 
epimorphism. It then follows that the strong epimorphisms are precisely
the {\em regular epimorphisms}; that is, the morphisms which are the 
coequalizer of some pair of maps. Our regular categories will always
be assumed to have finite limits. 

A full subcategory  is {\em reflective} when the inclusion has a left 
adjoint; this left adjoint is called the {\em reflection}. 

Throughout the paper, \S will be a category with finite limits; later
on we shall make further assumptions on \S, such as being 
regular; when we finally to come to our characterization of Grothendieck
quasitoposes, \S will be a presheaf topos.

Likewise, throughout the paper, \E will be a full reflective subcategory of \S.
We shall write $L$ for the reflection $\S\to\E$ and also sometimes for 
the induced endofunctor of \S, and we write $\ell\colon 1\to L$ for the
unit of the reflection. It is convenient to assume that the inclusion $\E\to\S$
is replete, in the sense that any object isomorphic to one in the image is itself
in the image. It is also convenient to assume that $\ell A\colon A\to LA$ is the 
identity whenever $A\in\E$. Neither assumption affects the results
of the paper. 

We shall say that the reflection has {\em monomorphic units} if each 
component $\ell X\colon X\to LX$ of the unit is a monomorphism, 
with an analogous meaning for {\em strongly epimorphic units}. 
When $L$ preserves finite limits it is said to be a {\em localization}. 

An object $A$ of a category \C is said to be {\em orthogonal} to a 
morphism
$f\colon X\to Y$ if each $a\colon X\to A$ factorizes uniquely through $f$. If 
instead each $a\colon X\to A$ factorizes in at most one way through $f$, the
object $A$ is said to be {\em separated with respect to $f$}, or 
{\em $f$-separated}. 
If \F is a class of morphisms, we say that $A$ is \F-orthogonal or \F-separated
if it is $f$-orthogonal or $f$-separated for each $f\in\F$. 

\section{Limit-preserving conditions for reflections}
\label{sect:conditions}

In this section we study various conditions on a reflection $L\colon\S\to\E$
weaker than being a localization. First observe that any reflective subcategory
is closed under limits, so the terminal object of \S lies in \E, and so $L$ always
preserves the terminal object. Thus preservation of finite limits is 
equivalent to preservation of pullbacks; our conditions all say that {\em certain}
pullbacks are preserved. 

\subsubsection*{\it Preservation of finite products.}
Since $L$ preserves the terminal object, preservation of finite products
amounts to preservation of binary products, or to preservation of 
pullbacks over the terminal object.

By a well-known result due to Brian Day \cite{Day-reflection}, if \S
is cartesian closed, then $L$ preserves finite products if and only 
if \E is an exponential ideal in \S; it then follows in particular that \E
is cartesian closed. 

\subsubsection*{\it Stable units.} 
For each object $B\in\E$, the reflection $L\colon\S\to\E$ induces
a reflection $L_B\colon\S/B\to\E/B$ onto the full subcategory 
$\E/B$ of $\S/B$. The original reflection $L$ is said to have
{\em stable units} when each $L_B$ preserves finite products,
or equivalently when $L$ preserves all pullbacks over objects 
of \E. Since the terminal object lies in \E, this implies in particular
that $L$ preserves finite products. 

If $\S$ is locally cartesian closed then, 
by the Day reflection theorem \cite{Day-reflection}
again, $L$ has stable units just
when each $\E/B$ is an exponential ideal in $\S/B$; it then follows
that \E is locally cartesian closed.

The name stable units was originally introduced in 
\cite{CassidyHebertKelly} for an apparently weaker condition, namely
that $L$ preserve each pullback of the form 
$$\xymatrix{
P \ar[r]^{q} \ar[d]_{p} & A \ar[d]^{u} \\
X \ar[r]_{\ell X} & LX }$$
but it was observed in 
\cite[Section~3.7]{CJKP-LocalizationStabilization} that these 
two conditions are in fact equivalent. Notice also that 
since $L\ell X$ is invertible, to say that $L$ preserves the 
pullback is equivalently to say that $L$ inverts $q$. 

\subsubsection*{\it Frobenius.} We say that $L$ satisfies the 
{\em Frobenius condition} when it preserves products 
of the form $X\x A$, with $A\in\E$. 

The condition is often given in the more general context
of an adjunction, not necessarily a reflection, between 
categories with finite products. In this case, the condition
is that the canonical map $$\phi:L(X\x A)\to LX\x A,$$ defined
using the comparison $L(X\x IA)\to LX\x LIA$ and the counit
$LIA\to A$, should be invertible. 

As is well-known, if \E and \S are both cartesian closed, then 
this condition is equivalent to the right adjoint $I:\E\to\S$ preserving
internal homs. As is perhaps less well-known, in our setting
of a reflection it is enough to assume
that \S is cartesian closed, and then the condition ensures that the
internal homs restrict to \E: see Proposition~\ref{prop:ccc} 
below. Thus if \S is cartesian closed, then $L$ satisfies the 
Frobenius conidition if and only if \E is closed in \S under
internal homs.

In fact, for any monadic adjunction satisfying the Frobenius
condition, internal homs may be lifted along the right adjoint.
More generally still, there is a version of the Frobenius condition
defined for monoidal categories in which the tensor product
is not required to be the product, and once again the internal
homs can be lifted along the right adjoint: 
see \cite[Proposition~3.5 and Theorem~3.6]{Hopf}. 

\subsubsection*{\it Semi-left-exact.}  
We say, following \cite{CassidyHebertKelly}, that $L$ 
is {\em semi-left-exact} if it preserves each pullback 
$$\xymatrix{
P \ar[r]^{q} \ar[d]_{p} & A \ar[d]^{u} \\
X \ar[r]_{\ell X} & LX }$$
with $A\in\E$. This is clearly implied by the stable units condition. 
By \cite[Theorem~4.3]{CassidyHebertKelly} it in fact implies,
and so is equivalent to, the apparently stronger condition that
$L$ preserve each pullback of the form 
$$\xymatrix{
P \ar[r]^{q} \ar[d]_{p} & A \ar[d]^{u} \\
X \ar[r]_{v} & B }$$
with $A$ and $B$ in \E. But this latter condition is in turn
equivalent to the condition that each $L_B\colon\S/B\to\E/B$ be
Frobenius. Thus we see that semi-left-exactness is in fact 
a ``localized'' version of the Frobenius condition. 
In particular, we may 
take $B=1$, and see that semi-left-exactness implies 
the Frobenius condition. 

Once again, if \S is locally cartesian closed, so that each 
$\S/B$ is cartesian closed, then $L$ is semi-left-exact 
just when each $\E/B$ is closed in $\S/B$ under internal 
homs. This implies that each $\E/B$ is cartesian closed,
and so that $\E$ is locally cartesian closed; see 
Lemma~\ref{lemma:lccc} below. 

\subsubsection*{\it Preservation of  monomorphisms.} 
Preservation of monomorphisms will also be an important 
condition in what follows. Once again it can be 
seen as preservation of certain pullbacks. Notice also that 
$L\colon \S\to\E$ satisfies this condition if and only if each 
$L_B\colon\S/B\to\E/B$  does so. 

\subsubsection*{\it Relationships between the conditions.}
We summarize in the diagram
$$\xymatrix{
\textnormal{stable units} \ar@{=>}[r]  \ar@{=>}[d] & 
\textnormal{semi-left-exact} \ar@{=>}[d]  & 
\textnormal{mono-preserving} \ar@{=>}[d]  \\
\textnormal{finite-product-preserving} \ar@{=>}[r] & 
\textnormal{Frobenius} & \textnormal{mono-preserving} 
}$$
the relationships found so far between these conditions. Each
condition in the top row amounts to requiring the condition 
below it to hold for all $L_B:\S/B\to\E/B$ with $B\in\E$. 

In Theorem~\ref{thm:equivalence} below, we shall see that if \S
is regular and $L$ preserves monomorphisms, then the stable
units condition is equivalent to the conjunction of the semi-left-exactness and finite-product-preserving conditions. In order to prove this,
we start by considering separately the case where the
components of the unit $\ell\colon 1\to L$ are monomorphisms 
and that where they are strong epimorphisms.

\begin{theorem}\label{thm:mono-case}
  If \S is finitely complete and the reflection $L\colon\S\to\E$ has 
monomorphic units, then the following are equivalent:
\begin{enumerate}[(i)]
\item $L$ preserves finite limits;
\item $L$ has stable units;
\item $L$ is semi-left-exact and preserves finite products.
\end{enumerate}
\end{theorem}

\proof
The implication $(i)\Rightarrow(ii)$ is trivial, while $(ii)\Rightarrow(iii)$ 
was observed above. Thus it remains to verify the implication 
$(iii)\Rightarrow(i)$. 

Suppose then that $L$ is semi-left-exact and preserves finite products.
We must show that it preserves equalizers.
Given $f,g\colon Y\rightrightarrows Z$ in \S, form the equalizer $e\colon X\to Y$ of $f$ and $g$,
and the equalizer $d\colon A\to LY$ of $Lf$ and $Lg$; of course $A\in\E$ since
\E is closed in \S under limits. There is a unique map $k\colon X\to A$ making
the diagram 
$$\xymatrix{
X \ar[r]^{e} \ar[d]_{k} & Y \ar@<1ex>[r]^{f} \ar@<-1ex>[r]_{g} \ar[d]_{\ell Y} &
Z \ar[d]^{\ell Z} \\
A \ar[r]_{d} & LY \ar@<1ex>[r]^{Lf} \ar@<-1ex>[r]_{Lg} & LZ }$$
commute. It follows easily from the fact that $\ell Z$ is a monomorphism
that the square on the left is a pullback. Since $A\in\E$, it follows by 
semi-left-exactness that $L$ inverts $k$, which is equivalently to 
say that $L$ preserves the equalizer of $f$ and $g$.
\endproof

\begin{theorem}\label{thm:epi-case}
  If \S is regular, and the reflection $L\colon\S\to\E$
has strongly epimorphic units, then the following are equivalent:
\begin{enumerate}[(i)]
\item $L$ preserves finite products and monomorphisms;
\item $L$ has stable units and preserves monomorphisms.
\end{enumerate}
\end{theorem}

\proof
Since having stable units always implies the preservation of finite products,
it suffices to show that $(i)$ implies $(ii)$. Suppose then that $L$ 
preserves finite products and monomorphisms.
Let 
$$\xymatrix{
P \ar[r]^q \ar[d]_p & Y \ar[d]^{u} \\ X \ar[r]_{\ell X} & LX }$$ 
be a pullback. Since $\ell X$ is a strong epimorphism,
so is its pullback $q$; but the left adjoint $L$ preserves strong epimorphisms
and so $Lq$ is also a strong epimorphism in \E. 

Since $L$ preserves finite products and monomorphisms, it preserves
jointly monomorphic pairs; thus $Lp$ and $Lq$ are jointly monomorphic. 
It follows that the canonical comparison from $LP$ to the pullback of 
$Lu$ and $L\ell X$ is a monomorphism, but this comparison is just $Lq$. Thus
$Lq$ is a strong epimorphism and a monomorphism, and so invertible.
This proves that $L$ has stable units. 
\endproof

Having understood separately the case of a reflection with monomorphic
units and that of one with strongly epimorphic units, we now combine
these to deal with the general situation.
The first step in this direction
is well-known; see \cite{CassidyHebertKelly} for example.

\begin{proposition}\label{prop:factorization}
Suppose that \S is regular.
If \R is the closure of \E in \S under subobjects, then the 
reflection of \S into \E factorizes as 
$$\xymatrix{
\E  \rtwocell^{\overline{L}}{`\bot} & \R \rtwocell^{L'}{`\bot} & \S 
}$$
where  $L'$ has strongly epimorphic units and
$\overline{L}$ has monomorphic units.
\end{proposition}

\proof
A straightforward argument shows that an object $X\in\S$ lies in \R if 
and only if the unit $\ell X\colon X\to LX$ is a monomorphism. Then 
the restriction $\overline{L}\colon \R\to\E$ of $L$ is clearly a reflection
of \R into \E.

Since \S is regular we may 
factorize $\ell\colon 1\to L$ as a strong epimorphism 
$\ell'\colon 1\to L'$ followed by a monomorphism $\kappa\colon L'\to L$. 
Since $L'X$ is a subobject of $LX$, it lies in \R. We claim that $\ell'X\colon X\to L'X$ is a reflection of $X$ into \R. Given an object $Y\in\R$, the unit
$\ell Y\colon Y\to LY$ is a monomorphism, and now if $f:X\to Y$ is
any morphism, then in the diagram
$$\xymatrix{
X \ar[r]^{\ell'X} \ar[dd]_{f} & L'X \ar[d]^{\kappa X} \\
& LX \ar[d]^{Lf} \\
Y \ar[r]_{\ell Y} & LY }$$
$\ell'X$ is a strong epimorphism and $\ell Y$ a monomorphism, so there
is a unique induced $g\colon L'X\to Y$ with $g.\ell'X=f$ and $\ell Y.g=Lf.\kappa X$. This gives the existence of a factorization of $f$ through $\ell'X$;
the uniqueness is automatic since $\ell'X$ is a (strong) epimorphism. 
\endproof

\begin{corollary}\label{cor:LfromRtoE}
Suppose that \S is regular. 
If $L\colon\S\to\E$ is semi-left-exact and preserves finite products and monomorphisms, then its 
restriction  $\overline{L}\colon\R\to\E$  to \R preserves
finite limits, while $L'\colon\S\to\R$ has stable units and preserves
monomorphisms.
\end{corollary}

\proof Since $L$ is semi-left-exact and preserves finite products and mono\-morphisms, the same is true of its restriction $\overline{L}$. 
Thus $\overline{L}$ preserves finite limits by 
Theorem~\ref{thm:mono-case} and the fact that
$\overline{L}$ has monomorphic units. 

As for $L'$, since it has strongly epimorphic units it will suffice,
by Theorem~\ref{thm:epi-case}, to show that 
it preserves finite products and monomorphisms.

First observe that if $m\colon X\to Y$ is a monomorphism in \S, then
we have a commutative diagram
$$\xymatrix{
X \ar[r]^{\ell'X} \ar[d]_{m} & L'X \ar[r]^{\kappa X} \ar[d]_{L'm} & LX \ar[d]^{Lm} \\
Y \ar[r]_{\ell'Y} & L'Y \ar[r]_{\kappa Y} & LY }$$
in which $Lm$ and $\kappa X$ are monomorphisms, and thus also $L'm$.
This proves that $L'$ preserves monomorphisms. 

On the other hand, for any objects $X,Y\in\S$, we have a commutative
diagram 
$$\xymatrix @C3pc {
X\x Y \ar[r]^-{\ell'(X\x Y)} \ar[dr]_{\ell'X\x\ell'Y} & 
L'(X\x Y) \ar[d]^{\pi'} \ar[r]^{\kappa(X\x Y)} & L(X\x Y) \ar[d]^{\pi} \\
& L'X\x L'Y \ar[r]_{\kappa X\x \kappa Y} & LX\x LY }$$
in which $\pi$ and $\pi'$ are the canonical comparison maps. Now
$\ell'X\x \ell'Y$ is a strong epimorphism, since in a regular category
these are closed under products, and $\kappa X\x \kappa Y$ is a
monomorphism, since in any category these are closed under products. 
Since $L$ preserves products, $\pi$ is invertible, and it now follows 
that $\pi'$ is also invertible. Thus $L'$ preserves finite products. 
\endproof

\begin{theorem}\label{thm:equivalence}
Let \S be regular, and $L\colon \S\to\E$ an 
arbitrary reflection onto a full subcategory, with unit $\ell\colon 1\to L$. 
Then the following are equivalent:
\begin{enumerate}[(i)]
\item $L$ is semi-left-exact, and preserves finite products and monomorphisms;
\item $L$ has stable units and preserves monomorphisms. 
\end{enumerate}
\end{theorem}

\proof
The non-trivial part is that ($i$) implies ($ii$). Suppose then that $L$ 
satisfies ($i$), and factorize $L$ as $\overline{L}L'$ as in 
Proposition~\ref{prop:factorization}. By Corollary~\ref{cor:LfromRtoE},
we know that $\overline{L}$ preserves finite limits, while $L'$ preserves pullbacks over objects of \R and monomorphisms, thus the composite $\overline{L}L'$ preserves pullbacks over objects of \R and monomorphisms, and so in particular has stable units and preserves
monomorphisms. 
\endproof

\begin{remark}\label{rmk:stable-units}
In fact we have shown that a reflection $L$ satisfying the equivalent
conditions in the theorem preserves all pullbacks over an object of \R;
that is, over a subobject of an object in the subcategory~\E.
\end{remark}

Having described the positive relationships between our various 
conditions, we now show the extent to which they are independent. 
We shall give three examples; in each case \S is a presheaf category.

\begin{example}\label{ex}
  Let $\two$ be the full subcategory of $\Set\x\Set$ consisting of the
objects $(0,0)$ and $(1,1)$. This is reflective, with the reflection sending
a pair $(X,Y)$ to $(0,0)$ if $X$ and $Y$ are both empty, and $(1,1)$ 
otherwise. It's easy to see that this is semi-left-exact and 
preserves monomorphisms,
but fails to preserve the product $(0,1)\x(1,0)=(0,0)$ since $L(0,1)=L(0,1)=1$
but $L(0,0)=0$.
\end{example}

\begin{example}
Consider \Set as the full reflective subcategory of \RGph consisting 
of the discrete reflexive graphs. This time the reflector sends a graph $G$
to its set of connected components $\pi_0 G$. This is well-known to 
preserve finite products. Furthermore, it preserves pullbacks over a 
discrete reflexive graph $X$, since $\Set/X\simeq\Set^X$ and
$\RGph/X\simeq \RGph^X$ and the induced $\pi_0/X\colon \RGph/X\to\Set/X$
is just $\pi^X_0\colon \RGph^X\to\Set^X$, which preserves finite products since
$\pi_0$ does so. Thus $\pi_0$ has stable units, and so is semi-left-exact
(and, as we have already seen, preserves finite products). 
But $\pi_0$ does not preserve
monomorphisms: any set $X$ gives rise both to a discrete reflexive graph
and the ``complete'' reflexive graph $KX$ with exactly one directed 
edge between each pair of vertices. The inclusion 
$X\to KX$ is a monomorphism,
but $\pi_0X$ is just $X$, while $\pi_0KX=1$. Thus $\pi_0$ does not 
preserve this monomorphism if $X$ has more than one vertex. 
\end{example}

\begin{example}\label{ex:main}
Let \RGph be the category of reflexive graphs, and \Preord the full reflective
subcategory of preorders. Since the reflection sends a graph $G$ to 
a preorder on the set of vertices of $G$, it clearly preserves monomorphisms. 
An easy calculation shows that for any reflexive graph
$G$ and preorder $X$, the internal hom $[G,X]$ in \RGph again lies in \Preord,
corresponding to the set of graph homomorphisms equipped with the 
pointwise preordering.  Thus \Preord is an exponential ideal
in \RGph, and so the reflection preserves finite products by \cite{Day-reflection}. On the other hand, by Lemma~\ref{lemma:lccc} below,
the reflection cannot be semi-left-exact
since \Preord is not locally cartesian closed.
To see this, consider the preorder $X=\{x,y,y',z\}$
with $x\le y$ and $y'\le z$, and the two maps $1\to X$ picking out
$y$ and $y'$. Their coequalizer is $\{x\le y\le z\}$, but this is not 
preserved by pulling back along the inclusion of $\{x\le z\}$ into $X$, 
so \Preord cannot be locally cartesian closed. 
\end{example}

Our characterization of Grothendieck quasitoposes, in Theorem~\ref{thm:main} below, involves three conditions on a reflection: that it be 
semi-left-exact, that it preserve finite products, and that it preserve
monomorphisms. By the three examples above, we see that none of these three conditions can be omitted.

\section{Quasitoposes}\label{sect:qtop}

In this section
we take a slight detour to study  conditions under which a reflective
subcategory  is a quasitopos. First of all, a reflective subcategory \E of 
\S has any limits or colimits which \S does, so of course we have:

\begin{proposition}
If $L\colon \S\to\E$ is any reflection, then \E has finite limits and finite
colimits if \S does so. 
\end{proposition}

To deal with the remaining parts of the quasitopos structure we
require some assumptions on the reflection. 

\begin{proposition}\label{prop:ccc}
If $L\colon \S\to\E$ is Frobenius 
then \E is cartesian closed if \S is so.
\end{proposition}

\proof
Suppose that $L$ satisfies the Frobenius condition.
We shall show that if $A,B\in\E$, then $[A,B]$ is also in \E. 


The composite 
$$\xymatrix{
L[A,B]\x A \ar[r]^{\phi^{-1}} & L\left([A,B]\x A\right) \ar[r]^-{L\ev} & LB \ar[r]^\epsilon & B }$$
induces a morphism $c\colon L[A,B]\to [A,B]$. If we can show that 
$c\ell\colon [A,B]\to[A,B]$ is the identity, then $c$ will make $[A,B]$ into an 
$L$-algebra and so $[A,B]$ will lie in \E. 

Now commutativity of 
$$\xymatrix{
[A,B]\x A \ar[r]^{\ell\x 1} \ar[drr]_{\ell} \ar[dr]_{\ev} & L[A,B]\x A \ar[dr]^{\phi^{-1}} \\
& B \ar[dr]^{\ell} \ar[ddr]_1 & L\left([A,B]\x A\right) \ar[d]^{L\ev} \\
&& LB \ar[d]^{\epsilon}  \\
&& B }$$
shows that $\ev(c\ell\x1)=\ev$ and so that $c\ell=1$. \endproof


\begin{lemma}\label{lemma:lccc}
If $L\colon \S\to\E$ is semi-left-exact, then \E is locally cartesian closed if \S is so.
\end{lemma}

\proof
For any object $A\in\E$, the reflection $L$ induces a reflection of 
$\S/A$ into $\E/A$, which is Frobenius. It follows by
Proposition~\ref{prop:ccc}  that $\E/A$
is cartesian closed. \endproof

\begin{remark}
As observed in the previous section, there are converses to the 
previous two results. If \S is cartesian closed, and \E is a full
reflective subcategory closed under exponentials, then the 
reflection is Frobenius. And if \S is locally cartesian closed,
and \E is a full reflective subcategory closed under exponentials
in the slice categories, then the reflection is semi-left-exact.
\end{remark}

We now turn to the existence of weak subobject classifiers.
For this, we  consider one further condition on our 
reflection $L$,  weaker than preservation of finite limits. 
We say, following \cite{CarboniMantovani}, that
$L$ is {\em quasi-lex} if, for each finite diagram $X\colon \D\to\S$, the 
canonical comparison map $L(\lim X)\to \lim(LX)$ in \E is both a 
monomorphism and an epimorphism. We may then 
say that $L$ ``quasi-preserves'' the limit.



The proof of the next result closely follows that of 
\cite[Theorem~1.3.4]{CarboniMantovani}, although the assumptions 
made here are rather weaker. When we speak of 
unions of regular subobjects, we mean unions of subobjects which happen to 
be regular: there is no suggestion that the union itself must be regular. 
We say that such a union is {\em effective} when it is constructed as the 
pushout over the intersection.

\begin{proposition}\label{prop:quasi-lex}
Let \S be finitely complete and regular, and suppose further that \S 
has (epi, regular mono) factorizations of monomorphisms, and effective
unions of regular subobjects; for example, \S could be a quasitopos.
If  the reflection $L\colon\S\to\E$ preserves finite products and monomorphisms then it is quasi-lex. 
\end{proposition}

\proof
We know that $L$ preserves finite products, thus it will suffice to show 
that it quasi-preserves equalizers. 

Consider an equalizer diagram 
$$\xymatrix{
X \ar[r]^{e} & Y \ar@<1ex>[r]^{f} \ar@<-1ex>[r]_{g} & Z }$$
in \S. Since $L$ preserves finite products, quasi-preservation
of this equalizer is equivalent to quasi-preservation of
the equalizer 
$$\xymatrix{
X \ar[r]^{e} & Y \ar@<1ex>[r]^-{\binom{f}{Y}} \ar@<-1ex>[r]_-{\binom{g}{Y}} & Z\x Y }$$
in which now the parallel pair has a common retraction, given by the 
projection $Z\x Y\to Y$. This implies that the exterior of the diagram
$$\xymatrix{
X \ar[r]^{e} \ar[d]_{e} & Y \ar[d]_{f'} \ar[ddr]^{\binom{f}{Y}} \\
Y\ar[r]^{g'} \ar[drr]_{\binom{g}{Y}} & Z' \ar[dr]^(0.2){m} \\
&& Z\x Y }$$
is a pullback. By effectiveness of unions, we can form the union 
of $\binom{f}{Y}$ and $\binom{g}{Y}$ by constructing the pushout 
square as in the interior of the diagram. Then the induced map
$m\colon Z'\to Z\x Y$ will be the union, and in particular is a 
monomorphism. Now apply the reflection $L$ to this last diagram, to 
get a diagram 
$$\xymatrix{
LX \ar[r]^{Le} \ar[d]_{Le} & LY \ar[d]_{Lf'} \ar[ddr]^{\binom{Lf}{LY}} \\
LY\ar[r]^{Lg'} \ar[drr]_{\binom{Lg}{LY}} & LZ' \ar[dr]^(0.2){Lm} \\
&& LZ\x LY }$$
in \S. The interior square is still a pushout, and $Lm$ and $Le$ are still  
monomorphisms. Factorize $Le$ as an epimorphism $k\colon LX\to A$
followed by a regular monomorphism $d\colon A\to LY$. Then $d$, like
any regular monomorphism, is the equalizer of its cokernel pair. 
Since $k$ is an epimorphism, $d$ and $dk=Le$ have the same 
cokernel pair, namely $Lf'$ and $Lg'$. Thus $d$ is the equalizer of
$Lf'$ and $Lg'$, and $k$ is the canonical comparison. It is an epimorphism
by construction, and a monomorphism by the standard cancellation
properties. Thus $L$ quasi-preserves the equalizer of $f'$ and $g'$,
and so also the equalizer of $\binom{f}{Y}$ and $\binom{g}{Y}$, and so
finally that of $f$ and $g$.
\endproof

\begin{remark}
In fact there is also a partial converse to the preceding result: 
if $L$ is quasi-lex and has strongly
epimorphic units, then it preserves finite products and monomorphisms; indeed any quasi-lex $L$ preserves monomorphisms:
see \cite{CarboniMantovani}. 
\end{remark}

\begin{lemma}
If $L$ is quasi-lex, then \E has a weak subobject classifier if \S does so.
\end{lemma}

\proof
Let $t\colon 1\to\Omega$ be the weak subobject classifier of \S.
Now $Lt\colon L1\to L\Omega$ is a strong (in fact split) subobject, so there
is a unique map $\chi\colon L\Omega\to\Omega$ such that 
$$\xymatrix{
L1 \ar[r]^{Lt} \ar[d] & L\Omega \ar[d]^{\chi} \\
1 \ar[r]_{t} & \Omega}$$
is a pullback. Form the equalizer 
$$\xymatrix{
\Omega' \ar[r]^{i} & \Omega \ar@<1ex>[r]^{\chi\ell} \ar@<-1ex>[r]_{1} & 
\Omega }$$
in \S. 

Observe that $\chi.\ell.\chi=\chi.L\chi.\ell L\Omega=\chi.L\chi.L\ell\Omega$, and so $\chi.\ell.\chi.Li=\chi.L\chi.L\ell\Omega.Li=\chi.Li$; thus
$\chi.Li$ factorizes as $i.\chi'$ for a unique $\chi'\colon L\Omega'\to\Omega'$.

Furthermore, $i.\chi'.\ell\Omega'=\chi.Li.\ell\Omega'=\chi.\ell.i=i$
and so $\chi'.\ell=1$. This proves that $\Omega'\in\E$. Furthermore
$\chi.\ell.t=\chi.Lt.\ell=t$ and so $t=it'$ for a unique $t'\colon 1\to\Omega'$.
We shall show that $t'\colon 1\to\Omega'$  is a weak subobject classifier for \E. 

Suppose then that $m\colon A\to B$ is a strong subobject in \E. The inclusion,
being a right adjoint, preserves strong subobjects, so there is a unique
$f\colon B\to\Omega$ for which the diagram
$$\xymatrix{
A \ar[r]^{m} \ar[d] &  B \ar[d]^{f} \\
1 \ar[r]_{t} & \Omega }$$
is a pullback. We shall show that $f$ factorizes as $f=if'$; it then follows 
that $f'\colon B\to\Omega'$ is the unique map in \E classifying $m$.

To do so, it will suffice to show that $\chi.\ell.f=f$, or equivalently
$\chi.Lf.\ell=f$. Now consider the diagram
$$\xymatrix{
A \ar[r]^m \ar[d]_{\ell} & B \ar[d]^{\ell} \\
LA \ar[r]^{Lm} \ar[d]_{L!} & LB \ar[d]^{Lf} \\
L1 \ar[r]^{Lt} \ar[d] & L\Omega \ar[d]^{\chi} \\
1 \ar[r]_{t} & \Omega }$$
in which the top square is a pullback since $\ell A$ and $\ell B$
are invertible, and the bottom square is a pullback, by definition
of $\chi$. Thus if the middle square is a pullback, then the composite
will be, and so $\chi.Lf.\ell$ must be the unique map $f$ classifying $m$.

Now we know that the comparison $x$ from $LA$ to the pullback $P$
of $Lf$ and $Lt$ is both an epimorphism and a monomorphism in \E.
But $Lm$, like $m$, is a strong monomorphism, 
and factorizes as $sx$, where $s$ is the pullback of $Lt$;
thus $x$ is a strong monomorphism and an epimorphism, and so invertible.
This completes the proof.
\endproof




Combining the main results of this section, we have:

\begin{theorem}
If the reflection $L\colon \S\to\E$ is semi-left-exact and quasi-lex,
then \E is a quasitopos if \S is one. 
\end{theorem}

\begin{corollary}\label{cor:quasitopos}
  If the reflection $L\colon \S\to\E$ is semi-left-exact and preserves finite products and 
monomorphisms, and so also if it has stable units and preserves
monomorphisms, then \E is a quasitopos if \S is one.
\end{corollary}

\proof
Combine the previous theorem with Theorem~\ref{thm:equivalence}
and Proposition~\ref{prop:quasi-lex}.
\endproof

\section{Subquasitoposes}\label{sect:monics}

As recalled in the introduction, a {\em subtopos} of a topos is a full
reflective subcategory for which the reflection preserves finite limits. 
These can be characterized in various ways, for example using 
Lawvere-Tierney topologies, or universal closure operators. By analogy
with this, we define a {\em subquasitopos} of a quasitopos \S to be 
a full reflective subcategory for which the reflection has stable
units and preserves monomorphisms. 
By Corollary~\ref{cor:quasitopos} we know that the subcategory 
will indeed be a quasitopos. In this section, we give a 
classification of subquasitoposes of \S using {\em proper
universal closure operators}. 

A closure operator $j$, on a category \C with finite limits, 
assigns to each subobject $A'\le A$ a subobject $j(A')\le A$ in 
such a way that $A'\le j(A')=j(j(A))$ and if $A_1\le A_2\le A$ then 
$j(A_1)\le j(A_2)\le A$. The closure operator is said to be {\em 
universal} if for each $f\colon B\to A$ and each $A'\le A$ 
we have $f^*(j(A'))=j(f^*(A'))$. It is said to be 
{\em proper}, especially in the case where \C is a quasitopos, if
$j(A')\le A$ is strong subobject whenever $A'\le A$ is one; of 
course this is automatic if \C is a topos, so that all subobjects
are strong.  If $j(A')\le A$ is a strong 
subobject for {\em all} subobjects 
$A'\le A$, then $j$ is said to be a {\em strict} 
universal closure operator.

Given a universal closure operator $j$ on \C, a subobject
$m\colon A'\to A$ is said to be {\em $j$-dense} if $j(A'\le A)=A$.
An object $X$ of \C is said to be a {\em $j$-sheaf} if it is orthogonal
to each $j$-dense monomorphism, and {\em $j$-separated} if it
is separated with respect to each $j$-dense morphism. 

Recall, for example from \cite[Theorem~A4.4.8]{elephant1}, that 
for a quasitopos \S there is a bijection between localizations of \S
and proper universal closure operators on \S. Explicitly,
the bijection associates to a proper universal closure operator
$j$ the subcategory $\Sh(\S,j)$ of $j$-sheaves; while for a localization
$L\colon\S\to\E$, the corresponding closure operator sends a 
subobject $A'\le A$ to the  pullback of $LA'\le LA$ along the 
unit $\ell\colon A\to LA$. Furthermore, by 
\cite[Theorem~A4.4.5]{elephant1}, if $j$ is strict then 
$\Sh(\S,j)$ is a topos. Conversely, if $j$ is a proper universal
closure operator for which $\Sh(\S,j)$ is a topos, then for 
any subobject $A'\le A$ in \S, the reflection $LA'\le LA$ is 
a subobject in a topos, hence a strong subobject; thus $j(A')\le A$
too is a strong suboject, and $j$ is strict. 

For any quasitopos \Q, there is a strict universal closure operator
sending a subobject $A'\le A$ to its strong closure $\overline{A'}\le A$, given by factorizing the inclusion $A'\to A$ as an epimorphism
$A'\to \overline{A'}$ followed by a strong monomorphism $\overline{A'}\to A$. An object of \Q is said to be {\em coarse} if it is a sheaf
for this closure operator, and we write $\Cs(\Q)$ for the full subcategory consisting of the coarse objects; this is a topos, and is reflective
in \Q via a finite-limit-preserving reflection
$\Q\to\Cs(\Q)$
which inverts precisely those monomorphisms which are
also epimorphisms; see \cite[A4.4]{elephant1}.

We now suppose that \S is a quasitopos, and $L\colon\S\to\E$
a reflection onto a subquasitopos.
As before, we write \R for the full subcategory of \S consisting
of those objects $X\in\S$ for which the unit $\ell\colon X\to LX$
is a monomorphism. 
We saw in Proposition~\ref{prop:factorization} that \R
is reflective in \S, and we saw in 
Corollary~\ref{cor:LfromRtoE} that this reflection has stable
units and preserves monomorphisms; thus by 
Corollary~\ref{cor:quasitopos} the category \R, 
like \E, is a quasitopos. 

\begin{proposition}\label{prop:K}
  There is a strict universal closure operator $k$ on \S whose 
sheaves are the coarse objects in \E and whose separated objects
are the objects of \R. The class \K of $k$-dense monomorphisms
consists of all those monomorphisms $m\colon X\to Y$ for which
the monomorphism $Lm\colon LX\to LY$ is also an epimorphism 
in \E. 
\end{proposition}

\proof
Write $H\colon\E\to\Cs(\E)$ for the reflection; by the remarks above 
it preserves finite limits. Recall from 
Proposition~\ref{prop:quasi-lex} that $L\colon\S\to\E$ is 
quasi-lex;  since $H$ preserves finite limits and  inverts the epimorphic monomorphisms, 
 the composite $HL\colon\S\to\Cs(\E)$
is a finite-limit-preserving reflection. It follows that there is a proper
universal closure operator $k$ on \S whose sheaves are the coarse
objects in \E. Since $\Cs(\E)$ is a topos, $k$ is strict. 

A monomorphism $m\colon X\to Y$ in \S is $k$-dense just when 
it is inverted
by $HL$; that is, just when the monomorphism $Lm$ is also 
an epimorphism. An object $A\in\E$ is certainly separated with
respect to such an $m$, since for any $a\colon X\to A$ the 
induced $La\colon LX\to A$ has at most one factorization through 
the epimorphism $Lm\colon LX\to LY$. Furthermore, the 
$m$-separated objects are closed under subobjects, so that every
object of \R is $k$-separated.

Conversely, suppose that $X$ is $k$-separated; that is, separated
with respect to each $k$-dense $m$. We must show that 
$\ell X\colon X\to LX$ is a monomorphism. 
Let $d,c\colon K\rightrightarrows X$ be 
the kernel pair of $\ell X$, and $\delta\colon X\to K$ the diagonal. If $X$ is 
$\delta$-separated, then since $d\delta=1=c\delta$, 
the two morphisms $d$ and $c$
must be equal, which is to say that $\ell X$ is a monomorphism. Thus 
it will suffice to show that $\delta$ is $k$-dense. 
Since $L$ preserves finite
products and monomorphisms, it also preserves jointly monomorphic 
pairs; thus $Ld$ and $Lc$ are, like $d$ and $c$, jointly monomorphic. 
On the other hand $L\ell X$ is invertible, and $L\ell X.Ld=L\ell X.Lc$,
and so $Ld=Lc$; thus in fact $Ld$ is monomorphic. But $L\delta$ is 
a section of $Ld$, and so both maps are invertible. In particular, since
$L\delta$ is invertible, $\delta$ is $k$-dense, and so $X\in\R$.
\endproof

We are now ready to prove our characterization of subquasitoposes.

\begin{theorem}\label{thm:k}
Subquasitoposes of a quasitopos \S are in bijection with pairs
$(h,k)$, where $k$ is a strict universal closure operator on \S,
and $h$ is a proper universal closure operator on $\Sep(\S,k)$
with the property that every $h$-dense subobject is also $k$-dense;
the subquasitopos corresponding to the pair $(h,k)$ is 
$\Sh(\Sep(\S,k),h)$. 
\end{theorem}

\proof 
If $k$ is a strict universal closure operator on \S, then the category
$\Sh(\S,k)$ of $k$-sheaves is reflective in \S via a finite-limit-preserving reflection $M$. The category $\Sep(\S,k)$ of $k$-separated
objects is also reflective, and we may obtain the reflection $M'$
by factorizing the unit $m\colon X\to MX$ of $M$ as a strong
epimorphism $m':X\to M'X$ followed by a monomorphism 
$\kappa\colon M'X\to MX$, exactly as in Proposition~\ref{prop:factorization}. By
Corollary~\ref{cor:LfromRtoE} we know that $M'$ has stable
units and preserves monomorphisms. Now
$\Sh(\Sep(\S,k),h)$ is reflective in $\Sep(\S,k)$ via a finite-limit-preserving reflection, and so the composite reflection 
$\S\to\Sh(\Sep(\S,k),h)$ has stable units and preserves monomorphisms.

Conversely, let $L\colon\S\to\E$ be a reflection onto a
subquasitopos.
As above, we define $k$ to be the strict universal closure operator
whose sheaves are the coarse objects in \E. 
By Corollary~\ref{cor:LfromRtoE}, we know that the restriction 
$\overline{L}\colon\R\to\E$ of $L$ to \R preserves finite limits,
and so corresponds to a proper universal closure operator 
$h$ on \R, whose category of sheaves is \E. Since every 
$k$-sheaf is an $h$-sheaf, every $h$-dense monomorphism is 
$k$-dense.

It remains to prove the uniqueness of the $h$ and $k$ giving rise
to $L\colon\S\to\E$ as in the first paragraph. 
We constructed
$M'$ above by factorizing $X\to MX$ as a strong epimorphism
$m'\colon X\to M'X$ followed by a monomorphism 
$\kappa\colon M'X\to MX$. Since
every $h$-dense monomorphism is $k$-dense, certainly every
$k$-sheaf is an $h$-sheaf. Thus $\kappa\colon M'X\to MX$
factorizes through $LX$ by some $\nu\colon M'X\to LX$,
necessarily monic, and now $\ell\colon X\to LX$ factorizes as
a strong epimorphsm $m'\colon X\to M'X$ followed by a 
monomorphism $\nu\colon M'X\to LX$. Thus $\Sep(\S,k)$
is uniquely determined by $L$. In general, there can be several
different proper universal closure operators with a given category
of separated objects, but by the discussion after
\cite[Theorem~A4.4.8]{elephant1}, there can be at most one 
{\em strict} universal closure operator with a given category of separated
objects. Thus $k$ is uniquely determined. Unlike the case of 
separated objects, a proper universal closure operator {\em is} 
uniquely
determined by its sheaves, and so $h$ is also uniquely determined.
\endproof

Observe that in our characterization the two proper universal
closure operators live on different categories. In the next 
section, we shall see that when \S is a presheaf topos, there is 
an alternative characterization in terms of two universal closure
operators on \S. In fact, even for a general quasitopos, we may give
a characterization purely in terms of structure existing in \S 
provided that we prepared to work with stable classes of monomorphisms rather than universal closure operators. 

As in Proposition~\ref{prop:K}, we let \K denote the class of 
monomorphisms $m\colon X\to Y$ in \S for which $Lm$ is an
epimorphism as well as a monomorphism; as there, these are the 
dense monomorphisms for a universal closure operator, and so in 
particular are stable under pullback. Now we let \J be 
the class of monomorphisms $m\colon X\to Y$ in \S, every pullback
of which is inverted by $L$. This is clearly the largest stable class of 
monomorphisms inverted by $L$. As we saw in Proposition~\ref{prop:factorization}, 
the unit $\ell X\colon X\to LX$ is a monomorphism for any $X\in\R$;
furthermore since $L$ has stable units, it preserves the pullback of $\ell X$
along any map, and so $L$ inverts not just $\ell X$ but also all of its pullbacks.
Thus $\ell X$ lies in \J for all $X\in\R$; more generally, since $L$
preserves all pullbacks over objects in \R by Remark~\ref{rmk:stable-units},
any monomorphism  $f\colon X\to Y$ in \R which is inverted by $L$ will lie 
in \J.

\begin{theorem}\label{thm:stable-systems}
Let \S be a quasitopos, and $L\colon\S\to\E$ a reflection 
onto a full subcategory. If $L$ has stable units and preserves 
monomorphisms, then 
\begin{enumerate}[(i)]
\item an
object $X$ of \S lies in \R just when it is \K-separated;
\item 
an object $X$ of \S lies
in \E just when it is \K-separated  and a \J-sheaf.
\end{enumerate}
\end{theorem}

\proof
We have already proved part ($i$) in Proposition~\ref{prop:K}.
For part ($ii$), first observe that if $A\in\E$ then $A$ is orthogonal 
 to all morphisms inverted by $L$, not just those in \J. 
It is of course also separated with respect to \K.

Conversely, if $A$ is \K-separated then it is in \R; but then 
$\ell A\colon A\to LA$
is in \J, and so if $A$ is a \J-sheaf then $\ell A$ must be invertible and so
$A\in\E$. \endproof

At the current level of generality, there seems no reason why \J 
need be the dense monomorphisms for a proper universal
closure operator on \S. In the following section we shall see that 
this will be so if \S is a presheaf topos.

\section{Grothendieck quasitoposes}\label{sect:Grothendieck}

In this final section we suppose that \S is a presheaf topos $[\C\op,\Set]$,
as well as the standing assumption that $L\colon\S\to\E$ is a 
reflection which preserves monomorphisms and has stable units. 
Recall that \J consists of the monomorphisms which are stably inverted
by the reflection $L$, and that \K consists of the monomorphisms $m$
for which $Lm$ is an epimorphism in \E as well as a monomorphism.
By Theorem~\ref{thm:k},  the class \K consists of the 
dense monomorphisms for a (proper) universal closure operator $k$
on \S; and by our new assumption that \S is a presheaf
topos, $k$ corresponds to a Grothendieck topology with 
the same sheaves and separated objects. Since at this stage we are 
really only interested in the sheaves and separated objects, we take
the liberty of using the same name $k$ for the topology as for the
universal closure operator. 

As for \J, since it is a stable system of monomorphisms, it 
can be seen as a coverage, in the sense of \cite{elephant1}, and
so generates a Grothendieck topology $j$ whose sheaves are
the objects orthogonal to \J.

\begin{theorem}\label{thm:main}
For a reflection  $L\colon [\C\op,\Set]\to\E$ onto a full subcategory of a
presheaf category, the following conditions are equivalent:
\begin{enumerate}[(i)]
\item The subcategory \E has the form $\Sep(k)\cap\Sh(j)$ for topologies
$j$ and $k$ on \C with $k$ containing $j$;
\item $L$ is semi-left-exact and preserves finite products and monomorphisms;
\item $L$ has stable units and preserves monomorphisms.
\end{enumerate}
\end{theorem}
\noindent
An \E as in the theorem is called a Grothendieck quasitopos; as we saw in the introduction, a category \E has this form for some
$\C$, $j$, and $k$ if and only if it is a locally presentable quasitopos
\cite{BorceuxPedicchio}. 

\proof
The equivalence of ($ii$) and ($iii$) was shown in Theorem~\ref{thm:equivalence}. The fact that these imply ($i$) now follows from 
Theorem~\ref{thm:stable-systems}. Thus it will suffice to suppose ($i$)
and show that ($iii$) follows. 

We have adjunctions 
$$\xymatrix{
\mathllap{\Sep(k)\cap}\Sh(j) \rtwocell^{L_2}{`{\bot}} & \Sh(j) \rtwocell^{L_1}{`{\bot}}  & 
[\C\op,\mathrlap{\Set]} }$$
and $L_1$ preserves all finite limits. It will clearly suffice to 
show that $L_2$ preserves monomorphisms as well as 
pullbacks over an object of $\Sep(k)\cap\Sh(j)$.

Now $\Sep(k)\cap\Sh(j)$ is just the category of separated objects
in the topos $\Sh(j)$ for a (Lawvere-Tierney) topology $k'$ in $\Sh(j)$.
Thus it will suffice to show that for a topos \S and a topology $k$,
the reflection $L\colon \S\to\Sep(k)$ preserves monomorphisms as 
well as pullbacks over separated 
objects. This is the  special case of (one direction of)  
Theorem~\ref{thm:k}, where $h$ is trivial.
\endproof

As we saw in the previous section, the topology $k$ 
can be recovered from $\Sep(k)\cap\Sh(j)$, since
$\Sh(k)$ is the topos of coarse objects in $\Sep(k)\cap\Sh(j)$, which
can be obtained by inverting all those morphisms in 
$\Sep(k)\cap\Sh(j)$ which are both 
monomorphisms and epimorphisms. Unlike the case of 
the (proper) universal closure operator $h$ of the previous section,
$j$ need not be uniquely determined, as we now explain.

There exist non-trivial topologies $k$ for which every separated object is a sheaf;
these were studied by Johnstone in \cite{PTJ-quintessential}. In this 
case, for {\em any} topology $j$ contained in $k$ we have 
$$\Sep(k)\cap\Sh(j)=\Sh(k)\cap\Sh(j)=\Sh(k),$$ 
where the last step holds since $\Sh(k)\subseteq\Sh(j)$.
In particular we could take $j$ to be either trivial or $k$ and obtain
the same subcategory $\Sh(k)$ as $\Sep(k)\cap\Sh(j)$. 

\begin{example}
For example, as explained in \cite[Example~A4.4.9]{elephant1}, 
we could take the category $\Set^M$ of $M$-sets, where $M$ is
the two-element monoid $M=\{1,e\}$, with $e^2=e$, or equivalently the 
category of sets equipped with an idempotent. Then \Set
can be seen as the full reflective subcategory of $M$-sets with
trivial action. The reflection $L\colon \Set^M\to\Set$ splits the idempotent;
this preserves all limits and so is certainly a localization. Since the
unit of the adjunction is epimorphic, every separated object for 
the induced topology $k$ is a sheaf. 
\end{example}

\begin{remark}\label{rmk:non-G}
  Theorem~\ref{thm:main} can be generalized to the case of a 
Grothen\-dieck topos \S in place of $[\C\op,\Set]$;
then $j$ and $k$ would be Lawvere-Tierney topologies on \S. 
It can further be generalized to the case where \S is a Grothendieck 
quasitopos, provided that we are willing to work with proper universal 
closure operators $j$ and $k$ rather than topologies. In either 
case, \E will still
be a quasitopos by Corollary~\ref{cor:quasitopos}, and is in fact a 
Grothendieck quasitopos. 
In the case of a quasitopos or topos \S which is not locally 
presentable, however, there seems no reason
why the objects orthogonal to \J should be the sheaves, either 
for a topology or a universal closure operator. 
\end{remark}



 \bibliographystyle{plain}

\end{document}